\documentclass[11pt]{article}

\usepackage{amsmath}
\usepackage{amssymb}
\usepackage{amsthm}
\usepackage{amscd}

\theoremstyle{plain}
\newtheorem{thm}{Theorem}
\newtheorem{prop}{Proposition}
\newtheorem{lem}{Lemma}

\theoremstyle{definition}

\theoremstyle{remark}
\newtheorem{rem}{Remark}

\newcommand{\slth}{\widehat{\mathfrak{sl}}_2}

\newcommand{\N}{{\mathbb N}}
\newcommand{\F}{{\mathcal F}}
\newcommand{\z}{{\zeta}}

\newcommand{\Z}{{\mathbb Z}}
\newcommand{\C}{{\mathbb C}}
\newcommand{\odd}{{\rm odd}}
\newcommand{\even}{{\rm even}}

\newcommand{\br}[1]{{\langle #1 \rangle}} 

\newcommand{\bea}{\begin{eqnarray}}
\newcommand{\ena}{\end{eqnarray}}
\newcommand{\be}{\begin{eqnarray*}}
\newcommand{\en}{\end{eqnarray*}}
\newcommand{\eg}{\end{gather*}}
\newcommand{\bg}{\begin{gather*}}
\newcommand{\bal}{\begin{align}}
\newcommand{\baln}{\begin{align*}}
\newcommand{\eal}{\end{align*}}

\newcommand{\non}{\nonumber}
\newcommand{\noi}{\noindent}
\newcommand{\lb}[1]{\label{#1}}

\begin{document}
\vspace*{21mm}
\begin{center}
{\Large\bf 
BRST-resolution for\\[2mm]
Principally Graded Wakimoto module of $\slth$
} 
\end{center}
\vspace{11mm}

\begin{center}
{\large Yuji Hara, Naihuan Jing, Kailash Misra}\\[4mm]
{\it Department of Mathematics, North Carolina State University, \\
Raleigh, NC 27695, USA}
\end{center}

\vspace{5mm}
\begin{abstract}
BRST-resolution is studied for the principally graded Wakimoto module 
of $\slth$ recently found in \cite{prcplW}. The submodule structure is completely determined and irreducible representations can be obtained as the zero-th cohomology group.
\end{abstract}

\renewcommand{\thefootnote}{\fnsymbol{footnote}}
\footnote[0]{e-mail:yhara@ms.u-tokyo.ac.jp}

\newpage
\section{Introduction}
Free field representation of infinite dimensional Lie algebras is a way to construct representations concretely. It also enables us to compute important quantities in representations. A well known example is solutions for the correlation functions in the Conformal Field Theory 
\cite{DF1}\cite{DF}.

Among many free field representations, the Feigin-Fuchs construction for the Virasoro algebra \cite{CT} and the Wakimoto representations for affine Lie algebras \cite{W}\cite{FeFr} are special in the sense that they can be used to realize representations of generic level. 
For generic highest weights, they are isomorphic to the Verma modules. When they are reducible, the irreducible representations are obtained by means of the BRST-cohomology \cite{Fel89}\cite{BeFe}\cite{FeFr}.
Complete description of their submodule
structures are also given with Jantzen filtration technique\cite{FeFu}.
 
Recently the Wakimoto representation in the {\em principal} gradation was  found for the first time for $\slth$ in \cite{prcplW}. 
In this paper, we study the submodule structure and relation to irreducible representations of this Wakimoto representation using the same technique used in the above studies. 
Since discussions and results are parallel to the preceding cases, we concentrate on the speciality of this Wakimoto representation in proofs.
Theorem \ref{thm:sub} is about submodule structure and the BRST-resolution is given in Theorem \ref{thm:cmplx}.

\section{Principally graded Wakimoto module of $\slth$}
In this section, we briefly recall $\slth$ and its Wakimoto representation in the principal gradation, for detail see \cite{prcplW}.

The affine Lie algebra $\slth$ in the principal gradation consists of the generators $\beta_n, x_n, c, \rho$ which satisfy the following relations:
\bea
&&c: \mbox{central},\non\\
&& \left[\rho,\beta_n  \right]=n \beta_{n},\qquad 
   \left[\rho,x_{n}  \right]=n x_{n},\non\\
&&\left[\beta_m,\beta_n\right]=mc \delta_{m+n,0},\non\\
&&\left[\beta_m,x_{n}\right]=2 x_{m+n},\non\\
&&\left[x_{ m},x_{n}\right]=
-2(-1)^m \beta_{m+n}+mc(-1)^m\delta_{m+n,0}.\non
\ena
We introduce currents $\beta(\zeta),x(\zeta)$ as
\begin{eqnarray*}
\beta(\zeta)=\sum_{n \in \Z} \beta_n \zeta^{-n},\qquad
x(\zeta)=\sum_{n \in \Z} x_{n} \zeta^{-n}.
\end{eqnarray*}
Let us introduce three kinds of bosonic currents
\begin{eqnarray*}
&&
\phi_0(\zeta)=-\sum_{n:\odd}
{\phi_{0,n}\over n} \zeta^{-n},\\
&&
\phi_1(\zeta)=-\sum_{n:\even\atop n\neq 0}
{\phi_{1,n}\over n} \zeta^{-n}+
\phi_{1,0} \ln \zeta + q,\\
&&
\phi_2(\zeta)=-\sum_{n:\odd}
{\phi_{2,n}\over n} \zeta^{-n},
\end{eqnarray*}
with commutation relations
\begin{eqnarray*}
&&[\phi_{0,m},\phi_{0,n}]=+(1-(-1)^m)2 km \delta_{m+n,0},\\
&&[\phi_{1,m},\phi_{1,n}]=+(1+(-1)^m)2 (k+2)m
\delta_{m+n,0},\\ 
&&[\phi_{2,m},\phi_{2,n}]=-(1-(-1)^m)2 km
\delta_{m+n,0},\\ &&[\phi_{1,0},q]=4(k+2).
\end{eqnarray*}
Fock spaces are constructed on the vacuum vector $\nu_0$ as
\be
&&
\phi_{i,n}\nu_0=0\quad\mbox{for}\;n\geq 0,\\
&&
\nu_{j}=\exp\left(\frac{j}{2(k+2)}q\right)\nu_0,\\
&&
\F_j=\C[\phi_{1,-2},\phi_{1,-4}\cdots\phi_{0,-1},\phi_{0,-3}\cdots\phi_{2,-1},\phi_{2,-3}\cdots]\nu_{j},
\en
We denote by $d:\F_{j,k}\rightarrow\F_{j,k}$ the grading operator
\be
[d,\phi_{i,n}]=-n\,\phi_{i,n},\quad 
[d,q]=\phi_{1,0},
\qquad d\nu_j=\frac{2j^2+k}{4(k+2)}\nu_j.
\en
Then the Wakimoto representation $(\pi_j,\F_j)$ is given by the following Theorem \cite{prcplW}.
\begin{thm}\lb{thm:prcplW}
Let $j,k$ be complex numbers with $k\neq 0,-2$. 
Then the following gives a level $k$ representation 
of $\slth$ on the Fock space $\F_{j}${\rm :}
\bea
&&
\pi_j\left(\beta(\zeta)\right)={1\over 2}\z\partial\phi_0(\zeta),
\lb{eq7}\\
&& 
\pi_j\left(x(\zeta)\right)={1\over 2}:
\left( \z\partial\phi_1(\zeta)+\z\partial\phi_2(\zeta) \right)
e^{{\phi_2(\zeta)\over k}+{\phi_0(\zeta)\over k}}:,
\lb{eq8}\\
&&
\pi_j(c)=k,\quad \pi_j(\rho)=-d. 
\lb{eq9}
\ena
\end{thm}

This representation is the highest weight in the sense that 
\be
&&\beta_n \nu_j=0,~x_n \nu_j=0\quad {\rm for}\; n>0,\\
&&x_0 \nu_j=j\nu_j. 
\en
The highest weight is 
\bea
\Lambda_j=\frac{k}{2}(\Lambda_1+\Lambda_0)+j(\Lambda_1-\Lambda_0)
-\frac{\{2j(j+k+2)+k(k+3)\}}{8(k+2)}\delta,
\ena
where $\Lambda_0,\Lambda_1$ are 
the fundamental weights of $\slth$ and $\delta=\alpha_0+\alpha_1$.  

\section{Main results}
Here and henceforth, we assume $k+2=p/p'$ where $p,p'$ are positive integers relatively prime to each other. The following notations are frequently used
\be
&&
j_{m,m'}=m-m'(k+2), \\
&&
\F_{m,m'}=\F_{j_{m,m'}},\\
&&
\Lambda_{m,m'}=\Lambda_{j_{m,m'}}.
\en

Among reducible $\F_j$'s, 
the case we study in this paper is when $j=j_{m\pm lp,\, m'+1/2}$
\bea
l\in\{0\}\cap{\mathbb N}, \;1\leq m\leq p-1, \;0\leq m'\leq p'-1.
\lb{eqn:degenerate}
\ena
For the Verma module $V(\Lambda_{m\pm lp,\, m'+1/2})$ having the same highest weight, this is the condition to be completely degenerate {\it i.e} the structure of Verma module is described by the following diagram \cite{BeFe}\cite{Kac}:
\bea
\setlength{\unitlength}{0.8mm}
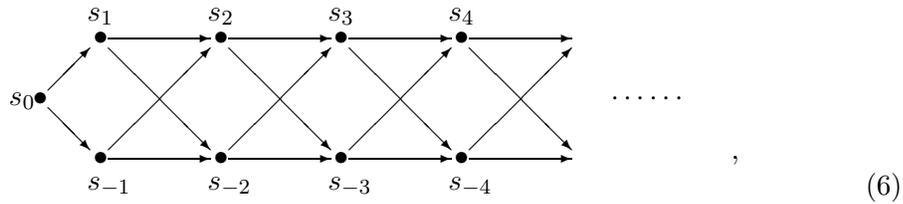
\begin{picture}(140,35)(0,4)
\put(9,19){$\bullet$}
\put(19,29){$\bullet$}
\put(19,9){$\bullet$}
\put(39,29){$\bullet$}
\put(39,9){$\bullet$}
\put(59,29){$\bullet$}
\put(59,9){$\bullet$}
\put(79,29){$\bullet$}
\put(79,9){$\bullet$}
\put(5,19){$s_0$}
\put(18,33){$s_{1}$}
\put(18,5){$s_{-1}$}
\put(38,33){$s_{2}$}
\put(38,5){$s_{-2}$}
\put(58,33){$s_{3}$}
\put(58,5){$s_{-3}$}
\put(78,33){$s_{4}$}
\put(78,5){$s_{-4}$}
\put(11.5,21.5){\vector(1,1){7}}
\put(11.5,18.5){\vector(1,-1){7}}
\put(21.5,30){\vector(1,0){17}}
\put(21.5,28.5){\vector(1,-1){17}}
\put(21.5,11.5){\vector(1,1){17}}
\put(21.5,10){\vector(1,0){17}}
\put(41.5,30){\vector(1,0){17}}
\put(41.5,28.5){\vector(1,-1){17}}
\put(41.5,11.5){\vector(1,1){17}}
\put(41.5,10){\vector(1,0){17}}
\put(61.5,30){\vector(1,0){17}}
\put(61.5,28.5){\vector(1,-1){17}}
\put(61.5,11.5){\vector(1,1){17}}
\put(61.5,10){\vector(1,0){17}}
\put(81.5,30){\vector(1,0){17}}
\put(81.5,28.5){\vector(1,-1){17}}
\put(81.5,11.5){\vector(1,1){17}}
\put(81.5,10){\vector(1,0){17}}
\put(105,19){$\cdots\cdots$}
\put(125,10){,}
\end{picture}
\lb{diag:Verma}
\ena
where $s_i$ denotes a singular vector. Throughout this paper, an arrow goes from one vector to another if and only if the second vector is in the $\slth$-submodule generated by the first one. Weights of $s_i$ are given as follows:
\be
&&
\mbox{weight}(s_{-2i+1})=\Lambda_{-m+(l+2i)p,\, m'+1/2},
\quad i\in\N,\\
&&
\mbox{weight}(s_{-2i})=\Lambda_{m-(l+2i)p,\, m'+1/2},
\quad i\in\N,\\
&&
\mbox{weight}(s_{2i})=\Lambda_{m+(l+2i)p,\, m'+1/2},
\quad i\in\N,\\
&&
\mbox{weight}(s_{2i+1})=\Lambda_{-m-(l+2i)p,\, m'+1/2},
\quad i\in\{0\}\cap\N.
\en

Using the screening current $S(\z)$ introduced in \cite{prcplW}, a screening charge $Q^n$ is defined as
\begin{align}
&
Q^n=
\int_{{\mathcal C}}\prod_{i=1}^n\frac{d\z_i^2}{\z_i^2}
S(\z_1)S(\z_2)\cdots S(\z_n),\\
&
S(\z)=
{1 \over  2}\zeta^{2 \over  k+2}:
\z\partial\phi_2(\zeta)
e^{-{\phi_1(\zeta)\over k+2}}:.
\end{align}
The contour ${\mathcal C}$ is the same one used for the screening charge in 
\cite{Fel89} where $v_i$ there should be replaced with $\z_i^2$.
The operator $Q^n$ is well defined on $\F_{m,m'}$ when 
\be
n\equiv m\mbox{ mod }p,
\en 
and an intertwiner of Fock modules:
\begin{align*}
&
Q^n:\F_{m,m'}\to\F_{m-2n,m'},\\
&
[Q^n,a]=0\mbox{ for }\forall a\in\slth.
\end{align*}

\begin{prop}\label{prop:sng-csng}
The Fock module $\F_{m\pm lp,\, m'+1/2}$ with \eqref{eqn:degenerate} has singular vectors 
$u_{i},\;i\in\mathbb{N}$ and cosingular vectors $w_{i},\;i\in\{0\}\cup\mathbb{N}$. 
Their weights coincide with the weights of $s_i\in V(\Lambda_{m,m'+1/2})$  as 
\begin{align}
&
{\rm weight}(u_{i})={\rm weight}(s_{-2i+1}),
\lb{eqn:wt-u}\\
&
{\rm weight}(w_{i})={\rm weight}(s_{2i+1}).
\lb{eqn:wt-w}
\end{align}
This is the complete list of all singular and cosingular vectors in the module.
\end{prop}

\begin{proof}
First we show the existence of cosingular vectors.
Let $\F_j^*$ be the restricted dual of $\F_j$ and  $\omega$  be the Chevalley involution:
\bal
&
\omega({\beta_n})=\beta_{-n},\non\\
&
\omega({x_n})=
\begin{cases}
& x_{-n},\quad n:\mbox{even},\non\\
& -x_{-n},\quad n:\mbox{odd}.\non
\end{cases}
\end{align}
We define the contragradient representation of $\slth$: $(\pi_j^*,\F_j^*)$ as
\begin{align*}
\br{\pi_j^*(a)v^*,v}=-\br{v^*,\pi_j\left(\omega(a)\right)v},\quad a\in\slth.
\lb{eqn:contra2}
\end{align*}
Throughout this proof $v,v^*$ denote arbitrary vectors of $\F_j,\F_j^*$. 
Then there is a one to one correspondence between singular (cosingular) vectors in $\F_j$ and cosingular (singular) vectors in $\F_j^*$ such that vectors in a pair have the same weight. 

Hence we study the singular vectors in $\F_j^*$.
The space $\F_j^*$ can be regarded as a Fock space as
\begin{gather*}
\F_j^*=\C[\,{}^t\phi_{1,-2},{}^t\phi_{1,-4}\cdots{}^t\phi_{0,-1},{}^t\phi_{0,-3}\cdots{}^t\phi_{2,-1},{}^t\phi_{2,-3}\cdots]\nu_j^*,\\
\br{\,{}^t\phi_{i,n}v^*,v}=\br{v^*,\phi_{i,n}v}.
\end{gather*}

Let $\sigma$ be the automorphism which changes $0-1$ index of the Chevalley generators of $\slth$:
\begin{align}
\sigma(\beta_n)=\beta_n,\quad \sigma(x_n)=-x_n.
\end{align}
Set $\pi'_{-j}=\pi_{j}\circ\sigma$, then $(\pi'_{-j},\F_j)$ is another Wakimoto representation with the highest weight $\Lambda_{-j}$.

Then we define an isomorphism $\varphi:(\pi'_{j},\F_{-j})\to(\pi_j^*,\F_j^*)$ by
\begin{align}
&
\varphi(\phi_{i,n})=(-1)^i\times{}^t\phi_{i,n},\\
&
\varphi\left(\pi'_{j}(a)\right)=\left(\pi_{j}^*(a)\right)\quad(a\in\slth).
\end{align}
The last equation comes from 
\begin{align*}
\br{\varphi\left(\pi'_{j}(a)\right)v^*,v}
=-\br{v^*,\pi_j\left(\omega(a)\right)v}
=\br{\pi_j^*(a)v^*,v}.
\end{align*}
With this $\varphi$ a singular vector in $(\pi'_{j},\F_{-j})$ can be mapped to a singular vector of the same weight in $(\pi_j^*,\F_j^*)$. 

Singular vectors in $(\pi'_{m\pm lp,m'+1/2},\F_{-m\mp lp,-m'-1/2})$ can be constructed explicitly with screening charges. Denoting the highest weight vector of $\F_{m,m'+1/2}$ by $\nu(m)$,
\begin{align}
Q^{ip-m}\;\nu\!\left((2i+l)p-m\right),\quad i\in\N,\lb{eqn:sng1}
\end{align}
is a singular vector in $\F_{m+lp,m'+1/2},\; l\in\N$.
 To verify that the vector \eqref{eqn:sng1} is non-vanishing, we checked that the matrix element
\begin{align*}
\br{\nu^*(m+lp),(\phi_{2,-2(m'+1/2)+2(l+i)p'})^{-m+ip}\:
Q^{ip-m}\: \nu\!\left((2i+l)p-m\right)},
\end{align*}
is non-zero by using a formula in Appendix A of \cite{DF}. Note that the quantity $m'+1/2$ makes the subscript $-2(m'+1/2)+2(l+i)p'$ be odd.

By the same argument, 
\begin{align}
Q^{(l+i)p-m}\;\nu\!\left((2i+l)p-m\right),\quad i\in\N,\lb{eqn:sng2}
\end{align}
is a singular vector
in $\F_{m-lp,m'+1/2},\; l\in\N$.
In $(\pi'_{m\pm lp,m'+1/2},\F_{-m\mp lp,-m'-1/2})$ vectors  
\eqref{eqn:sng1},\eqref{eqn:sng2} are singular vectors of the weight 
$\Lambda_{-m-(l+2i)p,\, m'+1/2},\quad i\in\{0\}\cap\N$ 
and this proves existence of $w_i$ and \eqref{eqn:wt-w}.

To conclude all cosingular vectors in  $(\pi_{m\pm lp,m'+1/2},\F_{m\pm lp,m'+1/2})$ are given by the argument above, we study the determinant
 $\det C(N,j)$. The matrix $C(N,j)$ connects the elements of
 $U(\slth)\nu_j$ and monomials of $\F_j$ at the degree $N$ \cite{KM}, {\it e.g.}
\begin{align}
\begin{pmatrix}
\beta_{-1}\nu_j \\
x_{-1}\nu_j 
\end{pmatrix}
=C(1,j)\times
\begin{pmatrix}
\phi_{0,-1}\nu_j \\
\phi_{2,-1}\nu_j. 
\end{pmatrix}
\lb{eqn:detC}
\end{align}
If we look at $\det C(N,j)$ as a polynomial of $j$, its zero of the first degree for the smallest $N$ corresponds to a cosingular vector of $(\pi_j,\F_j)$ of the degree $N$. 
\begin{lem}
Let $g(N)$ denote the number of monomials of the degree $N$ in $U(\slth)$, then 
\begin{align}
\det C(N,j)={\rm const.}\times
\prod_{
\begin{subarray}{c}
r\geq 1\\
s=1,3,5\cdots\\
rs\leq N
\end{subarray}
}
(j-j_{r,s/2})^{g(N-rs)}.
\end{align}
\end{lem}
\begin{proof}

The cosingular vectors $w_0$ of  $(\pi_{r,s/2},\F_{r,s/2})$
 give RHS as a divisor of $\det C(N,j)$.
The order of $\det C(N,j)$ is the number of $x_n$ ($n$:odd) in the LHS of \eqref{eqn:detC}, which is given by $\sum_{k=0}^N f(N,k)k$ where $f(N,k)$ is defined by the polynomial of indetermimates $x,y$ as
\begin{align*}
\prod_{a\geq 1}\frac{1}{(1-x^a)(1-yx^{2a-1})}=
\sum_{N\geq 0}\sum_{k=0}^N f(N,k)y^kx^N.
\end{align*}
And this order coincides with the order of RHS from the following identity:
\begin{align*}
\sum_{N\geq 0}\sum_{k=0}^N f(N,k)kx^N
&
=\sum_{a=1,3,5\cdots}\frac{x^a}{1-x^a}\prod_{b\geq 1}\frac{1}{(1-x^b)(1-x^{2b-1})}\\
&
=\sum_{N\geq 0}\sum_{
\begin{subarray}{c}
r\geq 1\\
s=1,3,5\cdots\\
rs\leq N
\end{subarray}
}
g(N-rs) x^N.
\end{align*}
\end{proof}

From this lemma, we can see that $\det C(N,j_{m\pm lp,m'+1/2})$ has zeros of the first degree only at $N=2(m+(l+i)p)(m'+1/2+ip')$ and $N=2(m+ip)(m'+1/2+(i+l)p')$ respectively. Since these zeros correspond to $w_i$, there is no cosingular vector other than $w_i$ .
This ends the proof for cosingular vectors $w_i$.


Using similar argument, we can also find all the cosingular vectors in 
$(\pi'_{m\pm lp},\F_{-m\mp lp,-m'-1/2})$
 and this proves that all singular vectors are given by $u_i$ and \eqref{eqn:wt-u}.
\end{proof}

\begin{rem}
The idea of using $\sigma$ is used in \cite{FeFr}.
\end{rem}

\begin{rem}
In the homogeneous gradation case, the automorphism $\sigma$ does not work well on the $\beta-\gamma$ system. The bosonization given in \cite{Frau} is much more suited to prove the similar proposition.
\end{rem}

\begin{thm}\label{thm:sub}
The submodule structure of $\F_{m\pm lp,\, m'+1/2}$ given by
\eqref{eqn:degenerate}
\begin{itemize}
\item[(1)]
The space of singular vectors in $\F_{m\pm lp,\, m'+1/2}$ is $\underset{1}{\overset{\infty}{\bigoplus}}\C u_{i}$. 
They generate a submodule ${\mathcal S}\F_{m\pm lp,\, m'+1/2}$ which is a direct sum  of irreducible highest weight modules.
The quotient $\F[1]=\F_{m\pm lp,\, m'+1/2}/{\mathcal S}\F_{m\pm lp,\, m'+1/2}$ contains infinitely many singular vectors $v_i\;(i\in\Z,\quad{\rm weight}(v_i)={\rm weight}(s_{2i}))$.
The space of singular vectors in $\F[1]$ is $\underset{-\infty}{\overset{\infty}{\bigoplus}}\C v_{i}$ and they generate a submodule ${\mathcal S}\F[1]$ which is also a direct sum  of irreducible highest weight modules.
Finally,  the space of singular vectors in $\F[2]=\F[1]/{\mathcal S}\F[1]$ is the projection of $\underset{0}{\overset{\infty}{\bigoplus}}\C w_{i}$ and $\F[2]$ is a direct sum  of irreducible highest weight modules.
\item[(2)]
Every submodule of $\F_{m\pm lp,\, m'+1/2}$ is generated by vectors $u_i$, $w_i$ and $v_i$.
\item[(3)]
The structure of the module $\F_{m\pm lp,\, m'+1/2}$ is described by this diagram
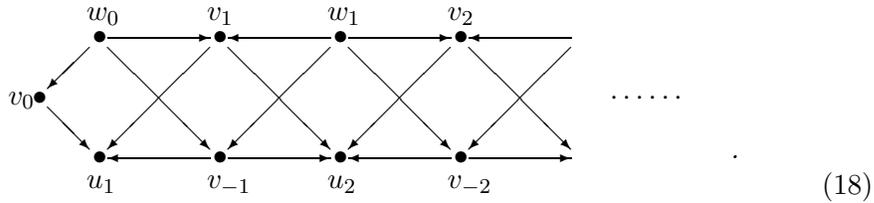
\bea
\setlength{\unitlength}{0.8mm}
\begin{picture}(140,35)(0,4)
\put(9,19){$\bullet$}
\put(19,29){$\bullet$}
\put(19,9){$\bullet$}
\put(39,29){$\bullet$}
\put(39,9){$\bullet$}
\put(59,29){$\bullet$}
\put(59,9){$\bullet$}
\put(79,29){$\bullet$}
\put(79,9){$\bullet$}
\put(5,19){$v_0$}
\put(18,33){$w_0$}
\put(18,5){$u_1$}
\put(38,33){$v_1$}
\put(38,5){$v_{-1}$}
\put(58,33){$w_1$}
\put(58,5){$u_2$}
\put(78,33){$v_2$}
\put(78,5){$v_{-2}$}
\put(18.5,28.5){\vector(-1,-1){7}}
\put(11.5,18.5){\vector(1,-1){7}}
\put(21.5,30){\vector(1,0){17}}
\put(21.5,28.5){\vector(1,-1){17}}
\put(38.5,28.5){\vector(-1,-1){17}}
\put(38.5,10){\vector(-1,0){17}}
\put(58.5,30){\vector(-1,0){17}}
\put(41.5,28.5){\vector(1,-1){17}}
\put(58.5,28.5){\vector(-1,-1){17}}
\put(41.5,10){\vector(1,0){17}}
\put(61.5,30){\vector(1,0){17}}
\put(61.5,28.5){\vector(1,-1){17}}
\put(78.5,28.5){\vector(-1,-1){17}}
\put(78.5,10){\vector(-1,0){17}}
\put(98.5,30){\vector(-1,0){17}}
\put(81.5,28.5){\vector(1,-1){17}}
\put(98.5,28.5){\vector(-1,-1){17}}
\put(81.5,10){\vector(1,0){17}}
\put(105,19){$\cdots\cdots$}
\put(125,10){.}
\end{picture}
\lb{diag:Wak}
\ena
\end{itemize}

\end{thm}
\begin{proof}
In the case of the homogeneous gradation \cite{BeFe}, 
to prove the same theorem the explicit form of Wakimoto construction is used only to show the existence of singular and cosingular vectors in Prop.4.5. Since their existence is proved in Prop.\ref{prop:sng-csng} in our case, the theorem can be proved using the same argument. 
\end{proof}

\begin{thm}\label{thm:cmplx}
\begin{itemize}
\item[]
\item[(1)]
If $1\leq m\leq p-1$ and $0\leq m'\leq p'-1$ then the following infinite
sequence
\[
\begin{CD}
\cdots @>{Q^{m}}>> \overset{-1}{\F_{2p-m,m'+\frac 12}} 
@>{Q^{p-m}}>> \overset{0}{\F_{m,m'+\frac 12}} 
@>{Q^{m}}>> \overset{1}{\F_{-m,m'+\frac 12}}
@>{Q^{p-m}}>> \cdots
\end{CD}
\]
is a complex $Q^{p-m}Q^m=Q^mQ^{p-m}=0$.
\item[(2)]
Its zero-th cohomology group is isomorphic to
the irreducible highest weight module $L(\Lambda_{m,m'+1/2})$
of the highest weight $\Lambda_{m,m'+\frac 12}${\rm :}
\be
\mbox{{\rm H}}^i=
\begin{cases}
& L(\Lambda_{m,m'+1/2}),\quad i=0,\\
& 0,\quad i\neq 0.
\end{cases}
\en
\end{itemize}
\end{thm}
\begin{proof}
This is proved by the same argument given in Section 4 of \cite{Fel89}. Note that as operators neither $Q^{p-m}Q^m$ nor $Q^mQ^{p-m}$ are zero but ${\rm Im}\,Q^m\subset {\rm Ker} \,Q^{p-m}$ and ${\rm Im} \,Q^{p-m} \subset {\rm Ker}\, Q^{m}$ hold on the Fock modules.

For (2), let $S':V(\Lambda_j)\to\F_j$ be the canonical homomorphism and modules $V_1,V_{-1},M_{v_0},M_{u_1}$ be submodules generated from $s_1,s_{-1},v_0,u_1$ in diagrams \eqref{diag:Verma}\eqref{diag:Wak}. Properties of the map $S'$ are known in the proof of Theorem \ref{thm:sub} with the Jantzen filtration and we can show
\begin{align*}
L(\Lambda_{m,m'+1/2})
\simeq S'\left(V(\Lambda_{m,m'+1/2})/(V_1+V_{-1})\right)
\simeq M_{v_0}/M_{u_1}
\simeq {\rm Ker}\, Q^{m}/{\rm Im} \,Q^{p-m}.
\end{align*}
\end{proof}

\noi {\Large{\it Acknowledgement}}\\
The first author would like to thank S. Odake for helpful discussions.



\end{document}